\documentclass[a4paper,11pt]{article}
\usepackage{rotating}
\usepackage{epsfig}
\usepackage{cite}
\usepackage{subfigure}
\usepackage{amsmath,amssymb,amscd}
\newtheorem{dfn}{Definition}
\newtheorem{thm}{Theorem}
\newtheorem{cor}{Corollary}

\newcommand{\mc}[1]{\mathcal{#1}}
\newcommand{\hsp}{\hspace{0.1in}}
\begin{document}
\begin{center}
{\huge Khovanov Homology and Conway Mutation}
\\[0.2in]
Stephan Wehrli, Universit\"at Basel, January 2003
\\[0.2in]
\end{center}
\begin{abstract} In this article, we present an easy example of mutant links with different
 Khovanov homology. The existence of such an example is important because it
 shows that Khovanov homology cannot be defined with a skein
rule similar to the skein relation for the Jones polynomial.
\end{abstract}

\section{Introduction}
In \cite{kh} M. Khovanov assigned to the diagram $D$ of an oriented link
$L$ a bigraded chain complex $\mc{C}^{*,*}(D)$, with a differential $d$ that
maps the chain group $\mc{C}^{i,j}(D)$ into $\mc{C}^{i+1,j}(D)$. 
He proved that the homotopy
equivalence class of graded chain complex $\mc{C}^{*,*}(D)$ only depends 
on the oriented link $L$. In
particular, the homology groups $\mc{H}^{i,j}(D)$ (considered up to isomorphism)
and the graded Poincar\'e polynomial 
$$ 
Kh(L)(t,q):=\sum_{i,j}t^i q^j
\operatorname{dim}_{\mathbb{Q}}(\mc{H}^{i,j}(D)\otimes\mathbb{Q}) 
\, \in \, \mathbb{Z}[t,t^{-1},q,q^{-1}]
$$ 
are link
invariants. The aim of this paper is to give an example of oriented mutant
links which are separated by the polynomial $Kh$ and to prove that consequently 
the invariant $Kh$ does not
satisfy a skein relation similar to the skein relation for the Jones polynomial.

\section{Skein equivalence} \label{skein}
In this section we briefly recall the definition of skein equivalence given
in \cite{ka}. A triple $(L_+,L_-,L_0)$ of oriented links is called a \emph{skein
triple}, if the oriented links $L_+$, $L_-$ and $L_0$ possess diagrams which are
mutually identical except in a small neighborhood, where they are 
respectively consistent with \begin{picture}(0,0)%
\includegraphics{l1.pstex}%
\end{picture}%
\setlength{\unitlength}{290sp}%
\begingroup\makeatletter\ifx\SetFigFont\undefined
\def\x#1#2#3#4#5#6#7\relax{\def\x{#1#2#3#4#5#6}}%
\expandafter\x\fmtname xxxxxx\relax \def\y{splain}%
\ifx\x\y   
\gdef\SetFigFont#1#2#3{%
  \ifnum #1<17\tiny\else \ifnum #1<20\small\else
  \ifnum #1<24\normalsize\else \ifnum #1<29\large\else
  \ifnum #1<34\Large\else \ifnum #1<41\LARGE\else
     \huge\fi\fi\fi\fi\fi\fi
  \csname #3\endcsname}%
\else
\gdef\SetFigFont#1#2#3{\begingroup
  \count@#1\relax \ifnum 25<\count@\count@25\fi
  \def\x{\endgroup\@setsize\SetFigFont{#2pt}}%
  \expandafter\x
    \csname \romannumeral\the\count@ pt\expandafter\endcsname
    \csname @\romannumeral\the\count@ pt\endcsname
  \csname #3\endcsname}%
\fi
\fi\endgroup
\begin{picture}(2464,2464)(-106,-1518)
\end{picture}
, \begin{picture}(0,0)%
\includegraphics{l2.pstex}%
\end{picture}%
\setlength{\unitlength}{290sp}%
\begingroup\makeatletter\ifx\SetFigFont\undefined
\def\x#1#2#3#4#5#6#7\relax{\def\x{#1#2#3#4#5#6}}%
\expandafter\x\fmtname xxxxxx\relax \def\y{splain}%
\ifx\x\y   
\gdef\SetFigFont#1#2#3{%
  \ifnum #1<17\tiny\else \ifnum #1<20\small\else
  \ifnum #1<24\normalsize\else \ifnum #1<29\large\else
  \ifnum #1<34\Large\else \ifnum #1<41\LARGE\else
     \huge\fi\fi\fi\fi\fi\fi
  \csname #3\endcsname}%
\else
\gdef\SetFigFont#1#2#3{\begingroup
  \count@#1\relax \ifnum 25<\count@\count@25\fi
  \def\x{\endgroup\@setsize\SetFigFont{#2pt}}%
  \expandafter\x
    \csname \romannumeral\the\count@ pt\expandafter\endcsname
    \csname @\romannumeral\the\count@ pt\endcsname
  \csname #3\endcsname}%
\fi
\fi\endgroup
\begin{picture}(2464,2464)(-106,-1518)
\end{picture}
 and \begin{picture}(0,0)%
\includegraphics{l0.pstex}%
\end{picture}%
\setlength{\unitlength}{290sp}%
\begingroup\makeatletter\ifx\SetFigFont\undefined
\def\x#1#2#3#4#5#6#7\relax{\def\x{#1#2#3#4#5#6}}%
\expandafter\x\fmtname xxxxxx\relax \def\y{splain}%
\ifx\x\y   
\gdef\SetFigFont#1#2#3{%
  \ifnum #1<17\tiny\else \ifnum #1<20\small\else
  \ifnum #1<24\normalsize\else \ifnum #1<29\large\else
  \ifnum #1<34\Large\else \ifnum #1<41\LARGE\else
     \huge\fi\fi\fi\fi\fi\fi
  \csname #3\endcsname}%
\else
\gdef\SetFigFont#1#2#3{\begingroup
  \count@#1\relax \ifnum 25<\count@\count@25\fi
  \def\x{\endgroup\@setsize\SetFigFont{#2pt}}%
  \expandafter\x
    \csname \romannumeral\the\count@ pt\expandafter\endcsname
    \csname @\romannumeral\the\count@ pt\endcsname
  \csname #3\endcsname}%
\fi
\fi\endgroup
\begin{picture}(2464,2464)(-106,-1518)
\end{picture}
.
\begin{dfn}
The \emph{skein equivalence} is the minimal (with respect to set-theoretical inclusion)
equivalence relation "$\sim$" on the set of oriented links such that
\begin{enumerate}
\item $L\sim L'$ when $L$ and $L'$ are isotopic,
\item $L_0\sim L'_0$ and $L_-\sim L'_-$ imply $L_+\sim L'_+$,
\item $L_0\sim L'_0$ and $L_+\sim L'_+$ imply $L_-\sim L'_-$,
\end{enumerate} for any two skein triples $(L_+,L_-,L_0)$ and $(L'_+,L'_-,L'_0)$.
\end{dfn}
It is easy to see that such a minimal relation as postulated in the definition actually
exists. The definition is motivated by the following: Assume we are given an
invariant $f$ of oriented links, such as the Jones polynomial, which takes
values in an arbitrary ring $R$ and satisfies a relation
$$ \alpha f(L_+)+\beta f(L_-)+\gamma f(L_0)=0,
$$
where $\alpha, \beta\in R^*$ and $\gamma\in R$. Then $f(L_+)$ is determined by
$f(L_0)$ and $f(L_-)$, and $f(L_-)$ is determined by $f(L_0)$ and $f(L_+)$. The
minimality of "$\sim$" implies:
\begin{thm} \label{tskein}
Let $L$ and $L'$ be skein equivalent. Then $f(L)=f(L')$.
\end{thm}

\section{Conway mutation} \label{conway}
\begin{figure}[h]\begin{center}\input{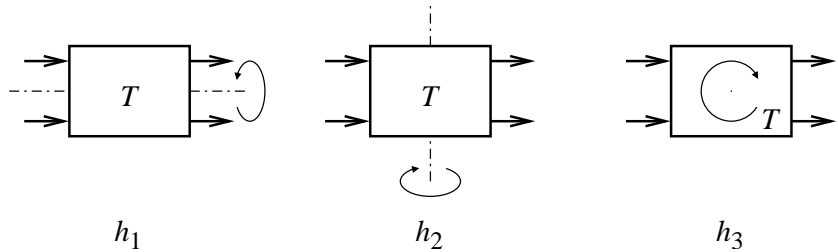}
\caption{The half-turns $h_1$, $h_2$ and $h_3$}\label{ftangleinv}\end{center}
\end{figure}

The mutation of links was originally defined in \cite{co}. We will use the
definition given in \cite{mu}. In Figure \ref{ftangleinv}, the rectangular boxes
represent an oriented 2-tangle $T$.
 Let $h_1$, $h_2$ and $h_3$ be the
half-turns about the indicated axes.
Define three involutions $\rho_1$, $\rho_2$ and $\rho_3$ on the
set of oriented 2-tangles by
$\rho_1T:=h_1(T)$, $\rho_2T:=-h_2(T)$ and $\rho_3T:=-h_3(T)$
(where $-h_2(T)$ and $-h_3(T)$ are the oriented 2-tangles
\begin{figure}[h]
\begin{center}\input{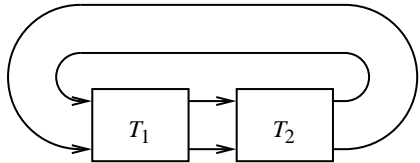}
\caption{The closure of the composition of $T_1$ and $T_2$}\label{ftanglecom} \end{center}
\end{figure}
obtained from $h_2(T)$ and $h_3(T)$ by reversing the orientations of all strings).
For two oriented 2-tangles $T_1$ and
$T_2$, denote by $T_1 T_2$ the composition of $T_1$ and $T_2$ and by
$(T_1 T_2)^\wedge$ the closure of $T_1 T_2$  (see Figure \ref{ftanglecom}).
\begin{dfn}
Two oriented links $L$ and $L'$ are called \emph{Conway mutants} if
there are two oriented 2-tangles $T_1$ and $T_2$ such that for an involution $\rho_i$
 $(i=1,2,3)$ the links $L$ and $L'$ are respectively isotopic to $(T_1 T_2)^\wedge$
 and $(T_1\rho_i T_2)^\wedge$.
\end{dfn}
\begin{thm} \label{tconway}
Let $L$ and $L'$ be Conway mutants. Then $L$ and $L'$ are skein equivalent.
\end{thm}
\emph{Proof.} The proof goes by induction on the number
$n$ of crossings of $T_2$. For $n\leq 1$, $T_2$ and $\rho_i T_2$ are isotopic,
whence $L\sim L'$. For $n>1$, modify a crossing of $T_2$ to obtain a skein
triple of tangles $(T_+,T_-,T_0)$ (with either $T_+=T_2$ or $T_-=T_2$, depending
on whether the crossing is positive or negative). Denote by $(L_+,L_-,L_0)$ and
$(L'_+,L'_-,L'_0)$ the skein triples corresponding to $(T_+,T_-,T_0)$ and
$(\rho_i T_+,\rho_i T_-,\rho_i T_0)$ respectively (i.e. $L_+=(T_1T_+)^\wedge$,
$L_-=(T_1T_-)^\wedge$ and so on). By induction, $L_0\sim L'_0$. Therefore, by
the definition of skein equivalence, $L_+\sim L'_+$ if and only if
$L_-\sim L'_-$. In other words, switching a crossing of $T_2$ does not affect
the truth or falsity of the assertion. Since $T_2$ can be untied by switching
crossings, we are back in the case $n\leq 1$. $\square$

\section{Mutant links with different Khovanov homology}
\begin{figure}
\input{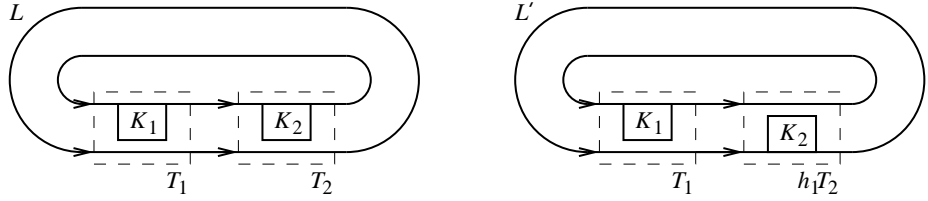} \caption{$L$ and $L'$ are Conway mutants}\label{fmutant}
\end{figure}

\begin{figure}[h]
\begin{center}\input{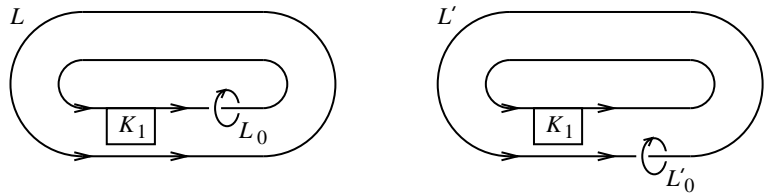}
\caption{$L$ and $L'$ for the case $n_2=2$} \label{fremark}
\end{center}
\end{figure}

Let $V(L)(q):=Kh(L)(-1,q)$ denote the graded Euler characteristic of $\mc{C}(D)$
and $W(L)(t):=Kh(L)(1,q)$ the ordinary (ungraded) Poincar\'e polynomial. As is
shown in \cite{kh}, $V$ is just an unnormalized version of the Jones polynomial.
By the results of sections \ref{skein} and \ref{conway}, the Jones polynomial is
invariant under Conway mutation. On the other hand, the following theorem gives
an example of mutant links which are separated by $W$.
\begin{thm} \label{tmain}
Let $K_i$ $(i=1,2)$ be a $(2,n_i)$ torus link, $n_i>2$. Then the oriented links
$$
L:=\begin{picture}(0,0)(0,500)%
\includegraphics{uk.pstex}%
\end{picture}%
\setlength{\unitlength}{290sp}%
\begingroup\makeatletter\ifx\SetFigFont\undefined
\def\x#1#2#3#4#5#6#7\relax{\def\x{#1#2#3#4#5#6}}%
\expandafter\x\fmtname xxxxxx\relax \def\y{splain}%
\ifx\x\y   
\gdef\SetFigFont#1#2#3{%
  \ifnum #1<17\tiny\else \ifnum #1<20\small\else
  \ifnum #1<24\normalsize\else \ifnum #1<29\large\else
  \ifnum #1<34\Large\else \ifnum #1<41\LARGE\else
     \huge\fi\fi\fi\fi\fi\fi
  \csname #3\endcsname}%
\else
\gdef\SetFigFont#1#2#3{\begingroup
  \count@#1\relax \ifnum 25<\count@\count@25\fi
  \def\x{\endgroup\@setsize\SetFigFont{#2pt}}%
  \expandafter\x
    \csname \romannumeral\the\count@ pt\expandafter\endcsname
    \csname @\romannumeral\the\count@ pt\endcsname
  \csname #3\endcsname}%
\fi
\fi\endgroup
\begin{picture}(2444,2754)(-96,-1508)
\end{picture}
\sqcup(K_1\sharp K_2)\\
\quad \text{and}\quad L':=K_1\sqcup K_2
$$
are Conway mutants with different $W$ polynomial. Here, $\begin{picture}(0,0)(0,500)%
\includegraphics{uk.pstex}%
\end{picture}%
\setlength{\unitlength}{290sp}%
\begingroup\makeatletter\ifx\SetFigFont\undefined
\def\x#1#2#3#4#5#6#7\relax{\def\x{#1#2#3#4#5#6}}%
\expandafter\x\fmtname xxxxxx\relax \def\y{splain}%
\ifx\x\y   
\gdef\SetFigFont#1#2#3{%
  \ifnum #1<17\tiny\else \ifnum #1<20\small\else
  \ifnum #1<24\normalsize\else \ifnum #1<29\large\else
  \ifnum #1<34\Large\else \ifnum #1<41\LARGE\else
     \huge\fi\fi\fi\fi\fi\fi
  \csname #3\endcsname}%
\else
\gdef\SetFigFont#1#2#3{\begingroup
  \count@#1\relax \ifnum 25<\count@\count@25\fi
  \def\x{\endgroup\@setsize\SetFigFont{#2pt}}%
  \expandafter\x
    \csname \romannumeral\the\count@ pt\expandafter\endcsname
    \csname @\romannumeral\the\count@ pt\endcsname
  \csname #3\endcsname}%
\fi
\fi\endgroup
\begin{picture}(2444,2754)(-96,-1508)
\end{picture}
$ denotes
the trivial knot and $K_1\sharp K_2$ is the connected sum of the oriented
links $K_1$ and $K_2$. Note that the connected sum is well-defined even if
$K_i$ has two components, because in this case the link $K_i$ is symmetric
in its components.
\end{thm}
\emph{Proof.} From Figure \ref{fmutant} it is apparent that $L$ and $L'$ are Conway
mutants. 
The Khovanov complex of the trivial knot is $$
\ldots\hsp\longrightarrow\hsp 0\hsp \longrightarrow\hsp
0\hsp\longrightarrow\hsp\mc{A}\hsp\longrightarrow\hsp 0\hsp\longrightarrow\hsp
0\hsp\longrightarrow\hsp\ldots $$ and $\operatorname{rank}(\mc{A})=2$, whence
$W(\begin{picture}(0,0)(0,500)%
\includegraphics{uk.pstex}%
\end{picture}%
\setlength{\unitlength}{290sp}%
\begingroup\makeatletter\ifx\SetFigFont\undefined
\def\x#1#2#3#4#5#6#7\relax{\def\x{#1#2#3#4#5#6}}%
\expandafter\x\fmtname xxxxxx\relax \def\y{splain}%
\ifx\x\y   
\gdef\SetFigFont#1#2#3{%
  \ifnum #1<17\tiny\else \ifnum #1<20\small\else
  \ifnum #1<24\normalsize\else \ifnum #1<29\large\else
  \ifnum #1<34\Large\else \ifnum #1<41\LARGE\else
     \huge\fi\fi\fi\fi\fi\fi
  \csname #3\endcsname}%
\else
\gdef\SetFigFont#1#2#3{\begingroup
  \count@#1\relax \ifnum 25<\count@\count@25\fi
  \def\x{\endgroup\@setsize\SetFigFont{#2pt}}%
  \expandafter\x
    \csname \romannumeral\the\count@ pt\expandafter\endcsname
    \csname @\romannumeral\the\count@ pt\endcsname
  \csname #3\endcsname}%
\fi
\fi\endgroup
\begin{picture}(2444,2754)(-96,-1508)
\end{picture}
)=2$. By [\citen{kh}, Proposition 33], $Kh$ is
multiplicative under disjoint union, and so $W(L)=2 W(K_1\sharp K_2)$. On the
other hand, by [\citen{kh}, Proposition 35],
$$
W(K_i)=2+t^{-2}+t^{-3}+\ldots+t^{-(n_i-1)}+t^{-n_i}
$$
if $n_i$ is odd and
$$
W(K_i)=2+t^{-2}+t^{-3}+\ldots+t^{-(n_i-1)}+2t^{-n_i}
$$
if $n_i$ is even. Since $n_i>2$, $W(K_i)$ is not divisible by 2. But
then $W(L')=W(K_1)W(K_2)$ is not divisible by 2 and hence $W(L')\neq W(L)$.
$\square$\\[11pt]
Theorems \ref{tskein}, \ref{tconway} and \ref{tmain} immediatly imply:
\begin{cor}
The $W$ polynomial does not satisfy a relation of the kind mentioned in
section \ref{skein}.
\end{cor}
\emph{Remark.} Theorem \ref{tmain} remains true if we also allow torus links
$K_i$ with $n_i<-2$ (this may be seen using [\citen{kh}, Corollary 11], which relates
the Khovanov homology of a link to the Khovanov homology of its mirror image).
The condition $|n_i|>2$ is necessary. In fact, if one of the
$|n_i|$ is $\leq 1$, then the corresponding torus link $K_i$ is trivial and hence $L$
and $L'$ are isotopic. If $n_2=2$, then $L$ and $L'$ look as is shown in Figure
\ref{fremark}. Note that both $L-L_0$ and $L'-L'_0$ are isotopic to the link
$\begin{picture}(0,0)(0,500)%
\includegraphics{uk.pstex}%
\end{picture}%
\setlength{\unitlength}{290sp}%
\begingroup\makeatletter\ifx\SetFigFont\undefined
\def\x#1#2#3#4#5#6#7\relax{\def\x{#1#2#3#4#5#6}}%
\expandafter\x\fmtname xxxxxx\relax \def\y{splain}%
\ifx\x\y   
\gdef\SetFigFont#1#2#3{%
  \ifnum #1<17\tiny\else \ifnum #1<20\small\else
  \ifnum #1<24\normalsize\else \ifnum #1<29\large\else
  \ifnum #1<34\Large\else \ifnum #1<41\LARGE\else
     \huge\fi\fi\fi\fi\fi\fi
  \csname #3\endcsname}%
\else
\gdef\SetFigFont#1#2#3{\begingroup
  \count@#1\relax \ifnum 25<\count@\count@25\fi
  \def\x{\endgroup\@setsize\SetFigFont{#2pt}}%
  \expandafter\x
    \csname \romannumeral\the\count@ pt\expandafter\endcsname
    \csname @\romannumeral\the\count@ pt\endcsname
  \csname #3\endcsname}%
\fi
\fi\endgroup
\begin{picture}(2444,2754)(-96,-1508)
\end{picture}
\sqcup K_1$. Using [\citen{kh}, Corollary 10],
one can show that both $\mc{H}^{i,j}(L)$ and $\mc{H}^{i,j}(L')$ are isomorphic to
$\mc{H}^{i+2,j+5}(\begin{picture}(0,0)(0,500)%
\includegraphics{uk.pstex}%
\end{picture}%
\setlength{\unitlength}{290sp}%
\begingroup\makeatletter\ifx\SetFigFont\undefined
\def\x#1#2#3#4#5#6#7\relax{\def\x{#1#2#3#4#5#6}}%
\expandafter\x\fmtname xxxxxx\relax \def\y{splain}%
\ifx\x\y   
\gdef\SetFigFont#1#2#3{%
  \ifnum #1<17\tiny\else \ifnum #1<20\small\else
  \ifnum #1<24\normalsize\else \ifnum #1<29\large\else
  \ifnum #1<34\Large\else \ifnum #1<41\LARGE\else
     \huge\fi\fi\fi\fi\fi\fi
  \csname #3\endcsname}%
\else
\gdef\SetFigFont#1#2#3{\begingroup
  \count@#1\relax \ifnum 25<\count@\count@25\fi
  \def\x{\endgroup\@setsize\SetFigFont{#2pt}}%
  \expandafter\x
    \csname \romannumeral\the\count@ pt\expandafter\endcsname
    \csname @\romannumeral\the\count@ pt\endcsname
  \csname #3\endcsname}%
\fi
\fi\endgroup
\begin{picture}(2444,2754)(-96,-1508)
\end{picture}
\sqcup K_1)\oplus \mc{H}^{i,j+1}(\begin{picture}(0,0)(0,500)%
\includegraphics{uk.pstex}%
\end{picture}%
\setlength{\unitlength}{290sp}%
\begingroup\makeatletter\ifx\SetFigFont\undefined
\def\x#1#2#3#4#5#6#7\relax{\def\x{#1#2#3#4#5#6}}%
\expandafter\x\fmtname xxxxxx\relax \def\y{splain}%
\ifx\x\y   
\gdef\SetFigFont#1#2#3{%
  \ifnum #1<17\tiny\else \ifnum #1<20\small\else
  \ifnum #1<24\normalsize\else \ifnum #1<29\large\else
  \ifnum #1<34\Large\else \ifnum #1<41\LARGE\else
     \huge\fi\fi\fi\fi\fi\fi
  \csname #3\endcsname}%
\else
\gdef\SetFigFont#1#2#3{\begingroup
  \count@#1\relax \ifnum 25<\count@\count@25\fi
  \def\x{\endgroup\@setsize\SetFigFont{#2pt}}%
  \expandafter\x
    \csname \romannumeral\the\count@ pt\expandafter\endcsname
    \csname @\romannumeral\the\count@ pt\endcsname
  \csname #3\endcsname}%
\fi
\fi\endgroup
\begin{picture}(2444,2754)(-96,-1508)
\end{picture}
\sqcup K_1)$.
The cases $n_2=-2$ and $n_1=\pm 2$ are similar.
\pagebreak
\section{Computer Calculations with KHOHO}
\begin{figure}[h]
\begin{center}\input{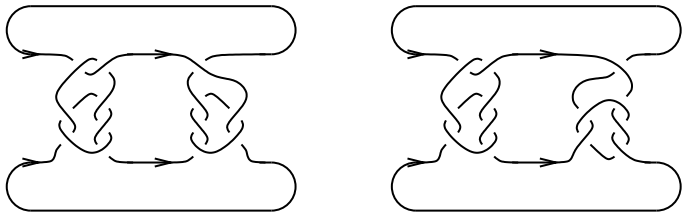}
\caption{The Kinoshita-Terasaka knot and the Conway knot}\label{fkinocon}\end{center}
\end{figure}

Tables \ref{t1} and \ref{t2} show the Khovanov homology of $L$ and
$L'$ for the case $n_1=n_2=3$. The tables where generated using A.
Shumakovitch's program KhoHo \cite{sh}. The entry in the $i$-th column and the
$j$-th row looks like $\hbox{$a$[$b$]}\over\lower 3pt\hbox{$c$}$, where $a$
is the rank of the homology group $\mc{H}^{i,j}$, $b$ the number of factors
$\mathbb{Z}/2 \mathbb{Z}$ in the decomposition of $\mc{H}^{i,j}$ into $p$-subgroups,
and $c$ the rank of the chain group $\mc{C}^{i,j}$. The numbers above the horizontal
arrows denote the ranks of the chain differentials.\\ In the examples,
only 2-torsion occurs. It has been conjectured by A. Shumakovitch that
this is actually true for arbitrary links. The reader may verify that not only
the dimensions but also the torsion parts of the $\mc{H}^{i,j}$ are different
for $L$ and $L'$.\\ The dimensions of the $\mc{C}^{i,j}$ agree because there is
a natural one-to-one correspondence between the Kauffman states of $L$ and $L'$
(which re-proves the fact that the Jones polynomial is invariant under Conway
mutation).\\ We do not know the answer to the following question:
\\ \noindent {\bf Question:} Does there exist a pair of mutant oriented knots with distinct
Khovanov homology?\\ According to the database of A. Shumakovitch, no such pair of knots with $13$
or less crossings exists. In particular, the Kinoshita-Terasaka knot and the Conway knot
(the knots depicted in Figure \ref{fkinocon}) are mutant knots with the same Khovanov homology
(see Table \ref{t3}).

\begin{table}[b]
\begin{center} 

\def\Cal#1{{\fam2#1}}

\let\TSp\thinspace
\let\NSp\negthinspace
\def\DSp{\thinspace\thinspace}
\def\QSp{\thinspace\thinspace\thinspace\thinspace}

\offinterlineskip
\def\gobble#1{}
\def\hline{\noalign{\hrule}}
\def\dblhline{height 0.16667em\gobble&&
&\vrule&
&\vrule&
&\vrule&
&\vrule&
&\vrule&
&\vrule&
\cr
\noalign{\hrule}}

\def\group#1#2#3#4{${\raise 1pt%
\hbox{#4#1\ifnum#2=0\else[#2]\fi}\over\lower 3pt\hbox{#4#3}}$}
\def\putdown#1{\smash{\vtop{\null\hbox{#1}}}}

\null\vfill
$$\vbox{\ialign{%
\vrule\TSp\vrule #\strut&\DSp\hfil #\DSp\vrule\TSp\vrule&\kern.75em
\TSp\hfil #\hfil\TSp&\hfil #\hfil&
\TSp\hfil #\hfil\TSp&\hfil #\hfil&
\TSp\hfil #\hfil\TSp&\hfil #\hfil&
\TSp\hfil #\hfil\TSp&\hfil #\hfil&
\TSp\hfil #\hfil\TSp&\hfil #\hfil&
\TSp\hfil #\hfil\TSp&\hfil #\hfil&
\TSp\hfil #\hfil\TSp\kern.75em\vrule\TSp\vrule\cr
\hline\dblhline
height 11pt depth 4pt&&
-6&\vrule&
-5&\vrule&
-4&\vrule&
-3&\vrule&
-2&\vrule&
-1&\vrule&
0\cr\hline\dblhline
height 0.2em depth 0.2em\gobble&&
&\vrule depth0pt&
&\vrule depth0pt&
&\vrule depth0pt&
&\vrule depth0pt&
&\vrule depth0pt&
&\vrule depth0pt&
\cr
height 0pt\gobble&&
&&
&&
&&
&&
&&
&&\cr
&-2&
&\smash{\vrule height 15pt}&
&\smash{\vrule height 15pt}&
&\smash{\vrule height 15pt}&
&\smash{\vrule height 15pt}&
&\smash{\vrule height 15pt}&
&\smash{\vrule height 15pt}&
\group{1}{0}{1}{\bf}\cr
height 0.3em\gobble&&
&\smash{\vrule height 10pt}&
&\smash{\vrule height 10pt}&
&\smash{\vrule height 10pt}&
&\smash{\vrule height 10pt}&
&\smash{\vrule height 10pt}&
&\smash{\vrule height 10pt}&
\cr
\hline
height 0.2em depth 0.2em\gobble&&
&\vrule depth0pt&
&\vrule depth0pt&
&\vrule depth0pt&
&\vrule depth0pt&
&\vrule depth0pt&
&\vrule depth0pt&
\cr
height 0pt\gobble&&
&&
&&
&&
&\QSp\putdown{2}\QSp&
&\QSp\putdown{4}\QSp&
&\QSp\putdown{2}\QSp&\cr
&-4&
&\smash{\vrule height 15pt}&
&\smash{\vrule height 15pt}&
&\smash{\vrule height 15pt}&
\group{0}{0}{2}{}&\rightarrowfill&
\group{0}{0}{6}{}&\rightarrowfill&
\group{0}{0}{6}{}&\rightarrowfill&
\group{2}{0}{4}{\bf}\cr
height 0.3em\gobble&&
&\smash{\vrule height 10pt}&
&\smash{\vrule height 10pt}&
&\smash{\vrule height 10pt}&
&\smash{\vrule height 10pt}&
&\smash{\vrule height 10pt}&
&\smash{\vrule height 10pt}&
\cr
\hline
height 0.2em depth 0.2em\gobble&&
&\vrule depth0pt&
&\vrule depth0pt&
&\vrule depth0pt&
&\vrule depth0pt&
&\vrule depth0pt&
&\vrule depth0pt&
\cr
height 0pt\gobble&&
&\QSp\putdown{1}\QSp&
&\QSp\putdown{5}\QSp&
&\DSp\putdown{10}\DSp&
&\DSp\putdown{18}\DSp&
&\DSp\putdown{13}\DSp&
&\QSp\putdown{5}\QSp&\cr
&-6&
\group{0}{0}{1}{}&\rightarrowfill&
\group{0}{0}{6}{}&\rightarrowfill&
\group{0}{0}{15}{}&\rightarrowfill&
\group{0}{0}{28}{}&\rightarrowfill&
\group{2}{0}{33}{\bf}&\rightarrowfill&
\group{0}{0}{18}{}&\rightarrowfill&
\group{1}{0}{6}{\bf}\cr
height 0.3em\gobble&&
&\smash{\vrule height 10pt}&
&\smash{\vrule height 10pt}&
&\smash{\vrule height 10pt}&
&\smash{\vrule height 10pt}&
&\smash{\vrule height 10pt}&
&\smash{\vrule height 10pt}&
\cr
\hline
height 0.2em depth 0.2em\gobble&&
&\vrule depth0pt&
&\vrule depth0pt&
&\vrule depth0pt&
&\vrule depth0pt&
&\vrule depth0pt&
&\vrule depth0pt&
\cr
height 0pt\gobble&&
&\QSp\putdown{6}\QSp&
&\DSp\putdown{24}\DSp&
&\DSp\putdown{36}\DSp&
&\DSp\putdown{38}\DSp&
&\DSp\putdown{14}\DSp&
&\QSp\putdown{4}\QSp&\cr
&-8&
\group{0}{0}{6}{}&\rightarrowfill&
\group{0}{0}{30}{}&\rightarrowfill&
\group{0}{0}{60}{}&\rightarrowfill&
\group{0}{0}{74}{}&\rightarrowfill&
\group{2}{2}{54}{\bf}&\rightarrowfill&
\group{0}{0}{18}{}&\rightarrowfill&
\group{0}{0}{4}{}\cr
height 0.3em\gobble&&
&\smash{\vrule height 10pt}&
&\smash{\vrule height 10pt}&
&\smash{\vrule height 10pt}&
&\smash{\vrule height 10pt}&
&\smash{\vrule height 10pt}&
&\smash{\vrule height 10pt}&
\cr
\hline
height 0.2em depth 0.2em\gobble&&
&\vrule depth0pt&
&\vrule depth0pt&
&\vrule depth0pt&
&\vrule depth0pt&
&\vrule depth0pt&
&\vrule depth0pt&
\cr
height 0pt\gobble&&
&\DSp\putdown{15}\DSp&
&\DSp\putdown{45}\DSp&
&\DSp\putdown{44}\DSp&
&\DSp\putdown{28}\DSp&
&\QSp\putdown{5}\QSp&
&\QSp\putdown{1}\QSp&\cr
&-10&
\group{0}{0}{15}{}&\rightarrowfill&
\group{0}{0}{60}{}&\rightarrowfill&
\group{1}{0}{90}{\bf}&\rightarrowfill&
\group{2}{0}{74}{\bf}&\rightarrowfill&
\group{0}{2}{33}{\bf}&\rightarrowfill&
\group{0}{0}{6}{}&\rightarrowfill&
\group{0}{0}{1}{}\cr
height 0.3em\gobble&&
&\smash{\vrule height 10pt}&
&\smash{\vrule height 10pt}&
&\smash{\vrule height 10pt}&
&\smash{\vrule height 10pt}&
&\smash{\vrule height 10pt}&
&\smash{\vrule height 10pt}&
\cr
\hline
height 0.2em depth 0.2em\gobble&&
&\vrule depth0pt&
&\vrule depth0pt&
&\vrule depth0pt&
&\vrule depth0pt&
&\vrule depth0pt&
&\vrule depth0pt&
\cr
height 0pt\gobble&&
&\DSp\putdown{20}\DSp&
&\DSp\putdown{39}\DSp&
&\DSp\putdown{20}\DSp&
&\QSp\putdown{6}\QSp&
&&
&&\cr
&-12&
\group{0}{0}{20}{}&\rightarrowfill&
\group{1}{0}{60}{\bf}&\rightarrowfill&
\group{1}{1}{60}{\bf}&\rightarrowfill&
\group{2}{0}{28}{\bf}&\rightarrowfill&
\group{0}{0}{6}{}&\smash{\vrule height 15pt}&
&\smash{\vrule height 15pt}&
\cr
height 0.3em\gobble&&
&\smash{\vrule height 10pt}&
&\smash{\vrule height 10pt}&
&\smash{\vrule height 10pt}&
&\smash{\vrule height 10pt}&
&\smash{\vrule height 10pt}&
&\smash{\vrule height 10pt}&
\cr
\hline
height 0.2em depth 0.2em\gobble&&
&\vrule depth0pt&
&\vrule depth0pt&
&\vrule depth0pt&
&\vrule depth0pt&
&\vrule depth0pt&
&\vrule depth0pt&
\cr
height 0pt\gobble&&
&\DSp\putdown{15}\DSp&
&\DSp\putdown{13}\DSp&
&\QSp\putdown{2}\QSp&
&&
&&
&&\cr
&-14&
\group{0}{0}{15}{}&\rightarrowfill&
\group{2}{1}{30}{\bf}&\rightarrowfill&
\group{0}{1}{15}{\bf}&\rightarrowfill&
\group{0}{0}{2}{}&\smash{\vrule height 15pt}&
&\smash{\vrule height 15pt}&
&\smash{\vrule height 15pt}&
\cr
height 0.3em\gobble&&
&\smash{\vrule height 10pt}&
&\smash{\vrule height 10pt}&
&\smash{\vrule height 10pt}&
&\smash{\vrule height 10pt}&
&\smash{\vrule height 10pt}&
&\smash{\vrule height 10pt}&
\cr
\hline
height 0.2em depth 0.2em\gobble&&
&\vrule depth0pt&
&\vrule depth0pt&
&\vrule depth0pt&
&\vrule depth0pt&
&\vrule depth0pt&
&\vrule depth0pt&
\cr
height 0pt\gobble&&
&\QSp\putdown{5}\QSp&
&&
&&
&&
&&
&&\cr
&-16&
\group{1}{0}{6}{\bf}&\rightarrowfill&
\group{1}{1}{6}{\bf}&\smash{\vrule height 15pt}&
&\smash{\vrule height 15pt}&
&\smash{\vrule height 15pt}&
&\smash{\vrule height 15pt}&
&\smash{\vrule height 15pt}&
\cr
height 0.3em\gobble&&
&\smash{\vrule height 10pt}&
&\smash{\vrule height 10pt}&
&\smash{\vrule height 10pt}&
&\smash{\vrule height 10pt}&
&\smash{\vrule height 10pt}&
&\smash{\vrule height 10pt}&
\cr
\hline
height 0.2em depth 0.2em\gobble&&
&\vrule depth0pt&
&\vrule depth0pt&
&\vrule depth0pt&
&\vrule depth0pt&
&\vrule depth0pt&
&\vrule depth0pt&
\cr
height 0pt\gobble&&
&&
&&
&&
&&
&&
&&\cr
&-18&
\group{1}{0}{1}{\bf}&\smash{\vrule height 15pt}&
&\smash{\vrule height 15pt}&
&\smash{\vrule height 15pt}&
&\smash{\vrule height 15pt}&
&\smash{\vrule height 15pt}&
&\smash{\vrule height 15pt}&
\cr
height 0.3em\gobble&&
&\smash{\vrule height 10pt}&
&\smash{\vrule height 10pt}&
&\smash{\vrule height 10pt}&
&\smash{\vrule height 10pt}&
&\smash{\vrule height 10pt}&
&\smash{\vrule height 10pt}&
\cr
\hline
\dblhline
}}$$


\caption{\vspace {4pt} Ranks of $\mc{H}^{i,j}$ and $\mc{C}^{i,j}$ and ranks
of the chain differentials for the disjoint union of the unknot and the
granny-knot}\label{t1}\end{center}
\end{table}

\begin{table}
\begin{center}

\def\Cal#1{{\fam2#1}}

\let\TSp\thinspace
\let\NSp\negthinspace
\def\DSp{\thinspace\thinspace}
\def\QSp{\thinspace\thinspace\thinspace\thinspace}

\offinterlineskip

\def\gobble#1{}
\def\hline{\noalign{\hrule}}
\def\dblhline{height 0.16667em\gobble&&
&\vrule&
&\vrule&
&\vrule&
&\vrule&
&\vrule&
&\vrule&
\cr
\noalign{\hrule}}

\def\group#1#2#3#4{${\raise 1pt%
\hbox{#4#1\ifnum#2=0\else[#2]\fi}\over\lower 3pt\hbox{#4#3}}$}
\def\putdown#1{\smash{\vtop{\null\hbox{#1}}}}

\null\vfill
$$\vbox{\ialign{%
\vrule\TSp\vrule #\strut&\DSp\hfil #\DSp\vrule\TSp\vrule&\kern.75em
\TSp\hfil #\hfil\TSp&\hfil #\hfil&
\TSp\hfil #\hfil\TSp&\hfil #\hfil&
\TSp\hfil #\hfil\TSp&\hfil #\hfil&
\TSp\hfil #\hfil\TSp&\hfil #\hfil&
\TSp\hfil #\hfil\TSp&\hfil #\hfil&
\TSp\hfil #\hfil\TSp&\hfil #\hfil&
\TSp\hfil #\hfil\TSp\kern.75em\vrule\TSp\vrule\cr
\hline\dblhline
height 11pt depth 4pt&&
-6&\vrule&
-5&\vrule&
-4&\vrule&
-3&\vrule&
-2&\vrule&
-1&\vrule&
0\cr\hline\dblhline
height 0.2em depth 0.2em\gobble&&
&\vrule depth0pt&
&\vrule depth0pt&
&\vrule depth0pt&
&\vrule depth0pt&
&\vrule depth0pt&
&\vrule depth0pt&
\cr
height 0pt\gobble&&
&&
&&
&&
&&
&&
&&\cr
&-2&
&\smash{\vrule height 15pt}&
&\smash{\vrule height 15pt}&
&\smash{\vrule height 15pt}&
&\smash{\vrule height 15pt}&
&\smash{\vrule height 15pt}&
&\smash{\vrule height 15pt}&
\group{1}{0}{1}{\bf}\cr
height 0.3em\gobble&&
&\smash{\vrule height 10pt}&
&\smash{\vrule height 10pt}&
&\smash{\vrule height 10pt}&
&\smash{\vrule height 10pt}&
&\smash{\vrule height 10pt}&
&\smash{\vrule height 10pt}&
\cr
\hline
height 0.2em depth 0.2em\gobble&&
&\vrule depth0pt&
&\vrule depth0pt&
&\vrule depth0pt&
&\vrule depth0pt&
&\vrule depth0pt&
&\vrule depth0pt&
\cr
height 0pt\gobble&&
&&
&&
&&
&\QSp\putdown{2}\QSp&
&\QSp\putdown{4}\QSp&
&\QSp\putdown{2}\QSp&\cr
&-4&
&\smash{\vrule height 15pt}&
&\smash{\vrule height 15pt}&
&\smash{\vrule height 15pt}&
\group{0}{0}{2}{}&\rightarrowfill&
\group{0}{0}{6}{}&\rightarrowfill&
\group{0}{0}{6}{}&\rightarrowfill&
\group{2}{0}{4}{\bf}\cr
height 0.3em\gobble&&
&\smash{\vrule height 10pt}&
&\smash{\vrule height 10pt}&
&\smash{\vrule height 10pt}&
&\smash{\vrule height 10pt}&
&\smash{\vrule height 10pt}&
&\smash{\vrule height 10pt}&
\cr
\hline
height 0.2em depth 0.2em\gobble&&
&\vrule depth0pt&
&\vrule depth0pt&
&\vrule depth0pt&
&\vrule depth0pt&
&\vrule depth0pt&
&\vrule depth0pt&
\cr
height 0pt\gobble&&
&\QSp\putdown{1}\QSp&
&\QSp\putdown{5}\QSp&
&\DSp\putdown{10}\DSp&
&\DSp\putdown{18}\DSp&
&\DSp\putdown{13}\DSp&
&\QSp\putdown{5}\QSp&\cr
&-6&
\group{0}{0}{1}{}&\rightarrowfill&
\group{0}{0}{6}{}&\rightarrowfill&
\group{0}{0}{15}{}&\rightarrowfill&
\group{0}{0}{28}{}&\rightarrowfill&
\group{2}{0}{33}{\bf}&\rightarrowfill&
\group{0}{0}{18}{}&\rightarrowfill&
\group{1}{0}{6}{\bf}\cr
height 0.3em\gobble&&
&\smash{\vrule height 10pt}&
&\smash{\vrule height 10pt}&
&\smash{\vrule height 10pt}&
&\smash{\vrule height 10pt}&
&\smash{\vrule height 10pt}&
&\smash{\vrule height 10pt}&
\cr
\hline
height 0.2em depth 0.2em\gobble&&
&\vrule depth0pt&
&\vrule depth0pt&
&\vrule depth0pt&
&\vrule depth0pt&
&\vrule depth0pt&
&\vrule depth0pt&
\cr
height 0pt\gobble&&
&\QSp\putdown{6}\QSp&
&\DSp\putdown{24}\DSp&
&\DSp\putdown{36}\DSp&
&\DSp\putdown{38}\DSp&
&\DSp\putdown{14}\DSp&
&\QSp\putdown{4}\QSp&\cr
&-8&
\group{0}{0}{6}{}&\rightarrowfill&
\group{0}{0}{30}{}&\rightarrowfill&
\group{0}{0}{60}{}&\rightarrowfill&
\group{0}{0}{74}{}&\rightarrowfill&
\group{2}{2}{54}{\bf}&\rightarrowfill&
\group{0}{0}{18}{}&\rightarrowfill&
\group{0}{0}{4}{}\cr
height 0.3em\gobble&&
&\smash{\vrule height 10pt}&
&\smash{\vrule height 10pt}&
&\smash{\vrule height 10pt}&
&\smash{\vrule height 10pt}&
&\smash{\vrule height 10pt}&
&\smash{\vrule height 10pt}&
\cr
\hline
height 0.2em depth 0.2em\gobble&&
&\vrule depth0pt&
&\vrule depth0pt&
&\vrule depth0pt&
&\vrule depth0pt&
&\vrule depth0pt&
&\vrule depth0pt&
\cr
height 0pt\gobble&&
&\DSp\putdown{15}\DSp&
&\DSp\putdown{45}\DSp&
&\DSp\putdown{44}\DSp&
&\DSp\putdown{28}\DSp&
&\QSp\putdown{5}\QSp&
&\QSp\putdown{1}\QSp&\cr
&-10&
\group{0}{0}{15}{}&\rightarrowfill&
\group{0}{0}{60}{}&\rightarrowfill&
\group{1}{0}{90}{\bf}&\rightarrowfill&
\group{2}{0}{74}{\bf}&\rightarrowfill&
\group{0}{2}{33}{\bf}&\rightarrowfill&
\group{0}{0}{6}{}&\rightarrowfill&
\group{0}{0}{1}{}\cr
height 0.3em\gobble&&
&\smash{\vrule height 10pt}&
&\smash{\vrule height 10pt}&
&\smash{\vrule height 10pt}&
&\smash{\vrule height 10pt}&
&\smash{\vrule height 10pt}&
&\smash{\vrule height 10pt}&
\cr
\hline
height 0.2em depth 0.2em\gobble&&
&\vrule depth0pt&
&\vrule depth0pt&
&\vrule depth0pt&
&\vrule depth0pt&
&\vrule depth0pt&
&\vrule depth0pt&
\cr
height 0pt\gobble&&
&\DSp\putdown{20}\DSp&
&\DSp\putdown{40}\DSp&
&\DSp\putdown{20}\DSp&
&\QSp\putdown{6}\QSp&
&&
&&\cr
&-12&
\group{0}{0}{20}{}&\rightarrowfill&
\group{0}{0}{60}{}&\rightarrowfill&
\group{0}{2}{60}{\bf}&\rightarrowfill&
\group{2}{0}{28}{\bf}&\rightarrowfill&
\group{0}{0}{6}{}&\smash{\vrule height 15pt}&
&\smash{\vrule height 15pt}&
\cr
height 0.3em\gobble&&
&\smash{\vrule height 10pt}&
&\smash{\vrule height 10pt}&
&\smash{\vrule height 10pt}&
&\smash{\vrule height 10pt}&
&\smash{\vrule height 10pt}&
&\smash{\vrule height 10pt}&
\cr
\hline
height 0.2em depth 0.2em\gobble&&
&\vrule depth0pt&
&\vrule depth0pt&
&\vrule depth0pt&
&\vrule depth0pt&
&\vrule depth0pt&
&\vrule depth0pt&
\cr
height 0pt\gobble&&
&\DSp\putdown{15}\DSp&
&\DSp\putdown{13}\DSp&
&\QSp\putdown{2}\QSp&
&&
&&
&&\cr
&-14&
\group{0}{0}{15}{}&\rightarrowfill&
\group{2}{1}{30}{\bf}&\rightarrowfill&
\group{0}{1}{15}{\bf}&\rightarrowfill&
\group{0}{0}{2}{}&\smash{\vrule height 15pt}&
&\smash{\vrule height 15pt}&
&\smash{\vrule height 15pt}&
\cr
height 0.3em\gobble&&
&\smash{\vrule height 10pt}&
&\smash{\vrule height 10pt}&
&\smash{\vrule height 10pt}&
&\smash{\vrule height 10pt}&
&\smash{\vrule height 10pt}&
&\smash{\vrule height 10pt}&
\cr
\hline
height 0.2em depth 0.2em\gobble&&
&\vrule depth0pt&
&\vrule depth0pt&
&\vrule depth0pt&
&\vrule depth0pt&
&\vrule depth0pt&
&\vrule depth0pt&
\cr
height 0pt\gobble&&
&\QSp\putdown{6}\QSp&
&&
&&
&&
&&
&&\cr
&-16&
\group{0}{0}{6}{}&\rightarrowfill&
\group{0}{2}{6}{\bf}&\smash{\vrule height 15pt}&
&\smash{\vrule height 15pt}&
&\smash{\vrule height 15pt}&
&\smash{\vrule height 15pt}&
&\smash{\vrule height 15pt}&
\cr
height 0.3em\gobble&&
&\smash{\vrule height 10pt}&
&\smash{\vrule height 10pt}&
&\smash{\vrule height 10pt}&
&\smash{\vrule height 10pt}&
&\smash{\vrule height 10pt}&
&\smash{\vrule height 10pt}&
\cr
\hline
height 0.2em depth 0.2em\gobble&&
&\vrule depth0pt&
&\vrule depth0pt&
&\vrule depth0pt&
&\vrule depth0pt&
&\vrule depth0pt&
&\vrule depth0pt&
\cr
height 0pt\gobble&&
&&
&&
&&
&&
&&
&&\cr
&-18&
\group{1}{0}{1}{\bf}&\smash{\vrule height 15pt}&
&\smash{\vrule height 15pt}&
&\smash{\vrule height 15pt}&
&\smash{\vrule height 15pt}&
&\smash{\vrule height 15pt}&
&\smash{\vrule height 15pt}&
\cr
height 0.3em\gobble&&
&\smash{\vrule height 10pt}&
&\smash{\vrule height 10pt}&
&\smash{\vrule height 10pt}&
&\smash{\vrule height 10pt}&
&\smash{\vrule height 10pt}&
&\smash{\vrule height 10pt}&
\cr
\hline
\dblhline
}}$$


\caption{\vspace{4pt} Ranks of $\mc{H}^{i,j}$ and $\mc{C}^{i,j}$ and ranks
of the chain differentials for the disjoint union of two trefoil knots}\label{t2}\end{center}
\end{table}

\begin{sidewaystable}
\begin{footnotesize}

\def\Cal#1{{\fam2#1}}

\let\TSp\thinspace
\let\NSp\negthinspace
\def\DSp{\thinspace\thinspace}
\def\QSp{\thinspace\thinspace\thinspace\thinspace}

\offinterlineskip

\def\gobble#1{}
\def\hline{\noalign{\hrule}}
\def\dblhline{height 0.16667em\gobble&&
&\vrule&
&\vrule&
&\vrule&
&\vrule&
&\vrule&
&\vrule&
&\vrule&
&\vrule&
&\vrule&
&\vrule&
&\vrule&
&\vrule&
&\vrule&
\cr
\noalign{\hrule}}

\def\group#1#2#3#4{${\raise 1pt%
\hbox{#4#1\ifnum#2=0\else[#2]\fi}\over\lower 3pt\hbox{#4#3}}$}
\def\putdown#1{\smash{\vtop{\null\hbox{#1}}}}

\null\vfill
$$\vbox{\ialign{%
\vrule\TSp\vrule #\strut&\DSp\hfil #\DSp\vrule\TSp\vrule&\kern.75em
\TSp\hfil #\hfil\TSp&\hfil #\hfil&
\TSp\hfil #\hfil\TSp&\hfil #\hfil&
\TSp\hfil #\hfil\TSp&\hfil #\hfil&
\TSp\hfil #\hfil\TSp&\hfil #\hfil&
\TSp\hfil #\hfil\TSp&\hfil #\hfil&
\TSp\hfil #\hfil\TSp&\hfil #\hfil&
\TSp\hfil #\hfil\TSp&\hfil #\hfil&
\TSp\hfil #\hfil\TSp&\hfil #\hfil&
\TSp\hfil #\hfil\TSp&\hfil #\hfil&
\TSp\hfil #\hfil\TSp&\hfil #\hfil&
\TSp\hfil #\hfil\TSp&\hfil #\hfil&
\TSp\hfil #\hfil\TSp&\hfil #\hfil&
\TSp\hfil #\hfil\TSp&\hfil #\hfil&
\TSp\hfil #\hfil\TSp\kern.75em\vrule\TSp\vrule\cr
\hline\dblhline
height 11pt depth 4pt&&
-7&\vrule&
-6&\vrule&
-5&\vrule&
-4&\vrule&
-3&\vrule&
-2&\vrule&
-1&\vrule&
0&\vrule&
1&\vrule&
2&\vrule&
3&\vrule&
4&\vrule&
5&\vrule&
6\cr\hline\dblhline
height 0.2em depth 0.2em\gobble&&
&\vrule depth0pt&
&\vrule depth0pt&
&\vrule depth0pt&
&\vrule depth0pt&
&\vrule depth0pt&
&\vrule depth0pt&
&\vrule depth0pt&
&\vrule depth0pt&
&\vrule depth0pt&
&\vrule depth0pt&
&\vrule depth0pt&
&\vrule depth0pt&
&\vrule depth0pt&
\cr
height 0pt\gobble&&
&&
&&
&&
&&
&&
&&
&&
&\QSp\putdown{1}\QSp&
&\QSp\putdown{5}\QSp&
&\DSp\putdown{10}\DSp&
&\DSp\putdown{10}\DSp&
&\QSp\putdown{6}\QSp&
&\QSp\putdown{1}\QSp&\cr
&9&
&\smash{\vrule height 15pt}&
&\smash{\vrule height 15pt}&
&\smash{\vrule height 15pt}&
&\smash{\vrule height 15pt}&
&\smash{\vrule height 15pt}&
&\smash{\vrule height 15pt}&
&\smash{\vrule height 15pt}&
\group{0}{0}{1}{}&\rightarrowfill&
\group{0}{0}{6}{}&\rightarrowfill&
\group{0}{0}{15}{}&\rightarrowfill&
\group{0}{0}{20}{}&\rightarrowfill&
\group{0}{0}{16}{}&\rightarrowfill&
\group{1}{0}{8}{\bf}&\rightarrowfill&
\group{0}{0}{1}{}\cr
height 0.3em\gobble&&
&\smash{\vrule height 10pt}&
&\smash{\vrule height 10pt}&
&\smash{\vrule height 10pt}&
&\smash{\vrule height 10pt}&
&\smash{\vrule height 10pt}&
&\smash{\vrule height 10pt}&
&\smash{\vrule height 10pt}&
&\smash{\vrule height 10pt}&
&\smash{\vrule height 10pt}&
&\smash{\vrule height 10pt}&
&\smash{\vrule height 10pt}&
&\smash{\vrule height 10pt}&
&\smash{\vrule height 10pt}&
\cr
\hline
height 0.2em depth 0.2em\gobble&&
&\vrule depth0pt&
&\vrule depth0pt&
&\vrule depth0pt&
&\vrule depth0pt&
&\vrule depth0pt&
&\vrule depth0pt&
&\vrule depth0pt&
&\vrule depth0pt&
&\vrule depth0pt&
&\vrule depth0pt&
&\vrule depth0pt&
&\vrule depth0pt&
&\vrule depth0pt&
\cr
height 0pt\gobble&&
&&
&&
&&
&&
&&
&\QSp\putdown{1}\QSp&
&\DSp\putdown{12}\DSp&
&\DSp\putdown{55}\DSp&
&\TSp\putdown{125}\TSp&
&\TSp\putdown{154}\TSp&
&\TSp\putdown{107}\TSp&
&\DSp\putdown{41}\DSp&
&\QSp\putdown{4}\QSp&\cr
&7&
&\smash{\vrule height 15pt}&
&\smash{\vrule height 15pt}&
&\smash{\vrule height 15pt}&
&\smash{\vrule height 15pt}&
&\smash{\vrule height 15pt}&
\group{0}{0}{1}{}&\rightarrowfill&
\group{0}{0}{13}{}&\rightarrowfill&
\group{0}{0}{67}{}&\rightarrowfill&
\group{0}{0}{180}{}&\rightarrowfill&
\group{0}{0}{279}{}&\rightarrowfill&
\group{0}{0}{261}{}&\rightarrowfill&
\group{1}{0}{149}{\bf}&\rightarrowfill&
\group{0}{1}{45}{\bf}&\rightarrowfill&
\group{0}{0}{4}{}\cr
height 0.3em\gobble&&
&\smash{\vrule height 10pt}&
&\smash{\vrule height 10pt}&
&\smash{\vrule height 10pt}&
&\smash{\vrule height 10pt}&
&\smash{\vrule height 10pt}&
&\smash{\vrule height 10pt}&
&\smash{\vrule height 10pt}&
&\smash{\vrule height 10pt}&
&\smash{\vrule height 10pt}&
&\smash{\vrule height 10pt}&
&\smash{\vrule height 10pt}&
&\smash{\vrule height 10pt}&
&\smash{\vrule height 10pt}&
\cr
\hline
height 0.2em depth 0.2em\gobble&&
&\vrule depth0pt&
&\vrule depth0pt&
&\vrule depth0pt&
&\vrule depth0pt&
&\vrule depth0pt&
&\vrule depth0pt&
&\vrule depth0pt&
&\vrule depth0pt&
&\vrule depth0pt&
&\vrule depth0pt&
&\vrule depth0pt&
&\vrule depth0pt&
&\vrule depth0pt&
\cr
height 0pt\gobble&&
&&
&&
&&
&&
&\QSp\putdown{5}\QSp&
&\DSp\putdown{55}\DSp&
&\TSp\putdown{260}\TSp&
&\TSp\putdown{658}\TSp&
&\TSp\putdown{951}\TSp&
&\TSp\putdown{792}\TSp&
&\TSp\putdown{371}\TSp&
&\DSp\putdown{89}\DSp&
&\QSp\putdown{6}\QSp&\cr
&5&
&\smash{\vrule height 15pt}&
&\smash{\vrule height 15pt}&
&\smash{\vrule height 15pt}&
&\smash{\vrule height 15pt}&
\group{0}{0}{5}{}&\rightarrowfill&
\group{0}{0}{60}{}&\rightarrowfill&
\group{0}{0}{315}{}&\rightarrowfill&
\group{0}{0}{918}{}&\rightarrowfill&
\group{0}{0}{1609}{}&\rightarrowfill&
\group{0}{0}{1743}{}&\rightarrowfill&
\group{1}{0}{1164}{\bf}&\rightarrowfill&
\group{1}{1}{461}{\bf}&\rightarrowfill&
\group{0}{0}{95}{}&\rightarrowfill&
\group{0}{0}{6}{}\cr
height 0.3em\gobble&&
&\smash{\vrule height 10pt}&
&\smash{\vrule height 10pt}&
&\smash{\vrule height 10pt}&
&\smash{\vrule height 10pt}&
&\smash{\vrule height 10pt}&
&\smash{\vrule height 10pt}&
&\smash{\vrule height 10pt}&
&\smash{\vrule height 10pt}&
&\smash{\vrule height 10pt}&
&\smash{\vrule height 10pt}&
&\smash{\vrule height 10pt}&
&\smash{\vrule height 10pt}&
&\smash{\vrule height 10pt}&
\cr
\hline
height 0.2em depth 0.2em\gobble&&
&\vrule depth0pt&
&\vrule depth0pt&
&\vrule depth0pt&
&\vrule depth0pt&
&\vrule depth0pt&
&\vrule depth0pt&
&\vrule depth0pt&
&\vrule depth0pt&
&\vrule depth0pt&
&\vrule depth0pt&
&\vrule depth0pt&
&\vrule depth0pt&
&\vrule depth0pt&
\cr
height 0pt\gobble&&
&&
&&
&&
&\DSp\putdown{10}\DSp&
&\TSp\putdown{126}\TSp&
&\TSp\putdown{663}\TSp&
&\TSp\putdown{1895}\TSp&
&\TSp\putdown{3159}\TSp&
&\TSp\putdown{3104}\TSp&
&\TSp\putdown{1767}\TSp&
&\TSp\putdown{565}\TSp&
&\DSp\putdown{91}\DSp&
&\QSp\putdown{4}\QSp&\cr
&3&
&\smash{\vrule height 15pt}&
&\smash{\vrule height 15pt}&
&\smash{\vrule height 15pt}&
\group{0}{0}{10}{}&\rightarrowfill&
\group{0}{0}{136}{}&\rightarrowfill&
\group{0}{0}{789}{}&\rightarrowfill&
\group{0}{0}{2558}{}&\rightarrowfill&
\group{0}{0}{5054}{}&\rightarrowfill&
\group{1}{0}{6264}{\bf}&\rightarrowfill&
\group{2}{0}{4873}{\bf}&\rightarrowfill&
\group{1}{1}{2333}{\bf}&\rightarrowfill&
\group{0}{0}{656}{}&\rightarrowfill&
\group{0}{0}{95}{}&\rightarrowfill&
\group{0}{0}{4}{}\cr
height 0.3em\gobble&&
&\smash{\vrule height 10pt}&
&\smash{\vrule height 10pt}&
&\smash{\vrule height 10pt}&
&\smash{\vrule height 10pt}&
&\smash{\vrule height 10pt}&
&\smash{\vrule height 10pt}&
&\smash{\vrule height 10pt}&
&\smash{\vrule height 10pt}&
&\smash{\vrule height 10pt}&
&\smash{\vrule height 10pt}&
&\smash{\vrule height 10pt}&
&\smash{\vrule height 10pt}&
&\smash{\vrule height 10pt}&
\cr
\hline
height 0.2em depth 0.2em\gobble&&
&\vrule depth0pt&
&\vrule depth0pt&
&\vrule depth0pt&
&\vrule depth0pt&
&\vrule depth0pt&
&\vrule depth0pt&
&\vrule depth0pt&
&\vrule depth0pt&
&\vrule depth0pt&
&\vrule depth0pt&
&\vrule depth0pt&
&\vrule depth0pt&
&\vrule depth0pt&
\cr
height 0pt\gobble&&
&&
&&
&\DSp\putdown{10}\DSp&
&\TSp\putdown{158}\TSp&
&\TSp\putdown{983}\TSp&
&\TSp\putdown{3241}\TSp&
&\TSp\putdown{6238}\TSp&
&\TSp\putdown{7166}\TSp&
&\TSp\putdown{4871}\TSp&
&\TSp\putdown{1916}\TSp&
&\TSp\putdown{417}\TSp&
&\DSp\putdown{44}\DSp&
&\QSp\putdown{1}\QSp&\cr
&1&
&\smash{\vrule height 15pt}&
&\smash{\vrule height 15pt}&
\group{0}{0}{10}{}&\rightarrowfill&
\group{0}{0}{168}{}&\rightarrowfill&
\group{0}{0}{1141}{}&\rightarrowfill&
\group{0}{0}{4224}{}&\rightarrowfill&
\group{0}{0}{9479}{}&\rightarrowfill&
\group{2}{0}{13406}{\bf}&\rightarrowfill&
\group{1}{1}{12038}{\bf}&\rightarrowfill&
\group{1}{2}{6788}{\bf}&\rightarrowfill&
\group{0}{0}{2333}{}&\rightarrowfill&
\group{0}{0}{461}{}&\rightarrowfill&
\group{0}{0}{45}{}&\rightarrowfill&
\group{0}{0}{1}{}\cr
height 0.3em\gobble&&
&\smash{\vrule height 10pt}&
&\smash{\vrule height 10pt}&
&\smash{\vrule height 10pt}&
&\smash{\vrule height 10pt}&
&\smash{\vrule height 10pt}&
&\smash{\vrule height 10pt}&
&\smash{\vrule height 10pt}&
&\smash{\vrule height 10pt}&
&\smash{\vrule height 10pt}&
&\smash{\vrule height 10pt}&
&\smash{\vrule height 10pt}&
&\smash{\vrule height 10pt}&
&\smash{\vrule height 10pt}&
\cr
\hline
height 0.2em depth 0.2em\gobble&&
&\vrule depth0pt&
&\vrule depth0pt&
&\vrule depth0pt&
&\vrule depth0pt&
&\vrule depth0pt&
&\vrule depth0pt&
&\vrule depth0pt&
&\vrule depth0pt&
&\vrule depth0pt&
&\vrule depth0pt&
&\vrule depth0pt&
&\vrule depth0pt&
&\vrule depth0pt&
\cr
height 0pt\gobble&&
&&
&\QSp\putdown{6}\QSp&
&\TSp\putdown{113}\TSp&
&\TSp\putdown{859}\TSp&
&\TSp\putdown{3375}\TSp&
&\TSp\putdown{7607}\TSp&
&\TSp\putdown{10219}\TSp&
&\TSp\putdown{8186}\TSp&
&\TSp\putdown{3850}\TSp&
&\TSp\putdown{1023}\TSp&
&\TSp\putdown{141}\TSp&
&\QSp\putdown{8}\QSp&
&&\cr
&-1&
&\smash{\vrule height 15pt}&
\group{0}{0}{6}{}&\rightarrowfill&
\group{0}{0}{119}{}&\rightarrowfill&
\group{0}{0}{972}{}&\rightarrowfill&
\group{0}{0}{4234}{}&\rightarrowfill&
\group{0}{0}{10982}{}&\rightarrowfill&
\group{1}{0}{17827}{\bf}&\rightarrowfill&
\group{3}{1}{18408}{\bf}&\rightarrowfill&
\group{2}{1}{12038}{\bf}&\rightarrowfill&
\group{0}{0}{4873}{}&\rightarrowfill&
\group{0}{0}{1164}{}&\rightarrowfill&
\group{0}{0}{149}{}&\rightarrowfill&
\group{0}{0}{8}{}&\smash{\vrule height 15pt}&
\cr
height 0.3em\gobble&&
&\smash{\vrule height 10pt}&
&\smash{\vrule height 10pt}&
&\smash{\vrule height 10pt}&
&\smash{\vrule height 10pt}&
&\smash{\vrule height 10pt}&
&\smash{\vrule height 10pt}&
&\smash{\vrule height 10pt}&
&\smash{\vrule height 10pt}&
&\smash{\vrule height 10pt}&
&\smash{\vrule height 10pt}&
&\smash{\vrule height 10pt}&
&\smash{\vrule height 10pt}&
&\smash{\vrule height 10pt}&
\cr
\hline
height 0.2em depth 0.2em\gobble&&
&\vrule depth0pt&
&\vrule depth0pt&
&\vrule depth0pt&
&\vrule depth0pt&
&\vrule depth0pt&
&\vrule depth0pt&
&\vrule depth0pt&
&\vrule depth0pt&
&\vrule depth0pt&
&\vrule depth0pt&
&\vrule depth0pt&
&\vrule depth0pt&
&\vrule depth0pt&
\cr
height 0pt\gobble&&
&\QSp\putdown{1}\QSp&
&\DSp\putdown{42}\DSp&
&\TSp\putdown{432}\TSp&
&\TSp\putdown{2128}\TSp&
&\TSp\putdown{5788}\TSp&
&\TSp\putdown{9186}\TSp&
&\TSp\putdown{8639}\TSp&
&\TSp\putdown{4766}\TSp&
&\TSp\putdown{1498}\TSp&
&\TSp\putdown{245}\TSp&
&\DSp\putdown{16}\DSp&
&&
&&\cr
&-3&
\group{0}{0}{1}{}&\rightarrowfill&
\group{0}{0}{43}{}&\rightarrowfill&
\group{0}{0}{474}{}&\rightarrowfill&
\group{0}{0}{2560}{}&\rightarrowfill&
\group{0}{0}{7916}{}&\rightarrowfill&
\group{2}{0}{14976}{\bf}&\rightarrowfill&
\group{2}{1}{17827}{\bf}&\rightarrowfill&
\group{1}{1}{13406}{\bf}&\rightarrowfill&
\group{0}{0}{6264}{}&\rightarrowfill&
\group{0}{0}{1743}{}&\rightarrowfill&
\group{0}{0}{261}{}&\rightarrowfill&
\group{0}{0}{16}{}&\smash{\vrule height 15pt}&
&\smash{\vrule height 15pt}&
\cr
height 0.3em\gobble&&
&\smash{\vrule height 10pt}&
&\smash{\vrule height 10pt}&
&\smash{\vrule height 10pt}&
&\smash{\vrule height 10pt}&
&\smash{\vrule height 10pt}&
&\smash{\vrule height 10pt}&
&\smash{\vrule height 10pt}&
&\smash{\vrule height 10pt}&
&\smash{\vrule height 10pt}&
&\smash{\vrule height 10pt}&
&\smash{\vrule height 10pt}&
&\smash{\vrule height 10pt}&
&\smash{\vrule height 10pt}&
\cr
\hline
height 0.2em depth 0.2em\gobble&&
&\vrule depth0pt&
&\vrule depth0pt&
&\vrule depth0pt&
&\vrule depth0pt&
&\vrule depth0pt&
&\vrule depth0pt&
&\vrule depth0pt&
&\vrule depth0pt&
&\vrule depth0pt&
&\vrule depth0pt&
&\vrule depth0pt&
&\vrule depth0pt&
&\vrule depth0pt&
\cr
height 0pt\gobble&&
&\QSp\putdown{5}\QSp&
&\TSp\putdown{113}\TSp&
&\TSp\putdown{784}\TSp&
&\TSp\putdown{2708}\TSp&
&\TSp\putdown{5207}\TSp&
&\TSp\putdown{5774}\TSp&
&\TSp\putdown{3704}\TSp&
&\TSp\putdown{1350}\TSp&
&\TSp\putdown{259}\TSp&
&\DSp\putdown{20}\DSp&
&&
&&
&&\cr
&-5&
\group{0}{0}{5}{}&\rightarrowfill&
\group{0}{0}{118}{}&\rightarrowfill&
\group{0}{0}{897}{}&\rightarrowfill&
\group{0}{0}{3492}{}&\rightarrowfill&
\group{1}{0}{7916}{\bf}&\rightarrowfill&
\group{1}{2}{10982}{\bf}&\rightarrowfill&
\group{1}{1}{9479}{\bf}&\rightarrowfill&
\group{0}{0}{5054}{}&\rightarrowfill&
\group{0}{0}{1609}{}&\rightarrowfill&
\group{0}{0}{279}{}&\rightarrowfill&
\group{0}{0}{20}{}&\smash{\vrule height 15pt}&
&\smash{\vrule height 15pt}&
&\smash{\vrule height 15pt}&
\cr
height 0.3em\gobble&&
&\smash{\vrule height 10pt}&
&\smash{\vrule height 10pt}&
&\smash{\vrule height 10pt}&
&\smash{\vrule height 10pt}&
&\smash{\vrule height 10pt}&
&\smash{\vrule height 10pt}&
&\smash{\vrule height 10pt}&
&\smash{\vrule height 10pt}&
&\smash{\vrule height 10pt}&
&\smash{\vrule height 10pt}&
&\smash{\vrule height 10pt}&
&\smash{\vrule height 10pt}&
&\smash{\vrule height 10pt}&
\cr
\hline
height 0.2em depth 0.2em\gobble&&
&\vrule depth0pt&
&\vrule depth0pt&
&\vrule depth0pt&
&\vrule depth0pt&
&\vrule depth0pt&
&\vrule depth0pt&
&\vrule depth0pt&
&\vrule depth0pt&
&\vrule depth0pt&
&\vrule depth0pt&
&\vrule depth0pt&
&\vrule depth0pt&
&\vrule depth0pt&
\cr
height 0pt\gobble&&
&\DSp\putdown{10}\DSp&
&\TSp\putdown{152}\TSp&
&\TSp\putdown{745}\TSp&
&\TSp\putdown{1814}\TSp&
&\TSp\putdown{2418}\TSp&
&\TSp\putdown{1805}\TSp&
&\TSp\putdown{753}\TSp&
&\TSp\putdown{165}\TSp&
&\DSp\putdown{15}\DSp&
&&
&&
&&
&&\cr
&-7&
\group{0}{0}{10}{}&\rightarrowfill&
\group{0}{0}{162}{}&\rightarrowfill&
\group{0}{0}{897}{}&\rightarrowfill&
\group{1}{0}{2560}{\bf}&\rightarrowfill&
\group{2}{1}{4234}{\bf}&\rightarrowfill&
\group{1}{0}{4224}{\bf}&\rightarrowfill&
\group{0}{0}{2558}{}&\rightarrowfill&
\group{0}{0}{918}{}&\rightarrowfill&
\group{0}{0}{180}{}&\rightarrowfill&
\group{0}{0}{15}{}&\smash{\vrule height 15pt}&
&\smash{\vrule height 15pt}&
&\smash{\vrule height 15pt}&
&\smash{\vrule height 15pt}&
\cr
height 0.3em\gobble&&
&\smash{\vrule height 10pt}&
&\smash{\vrule height 10pt}&
&\smash{\vrule height 10pt}&
&\smash{\vrule height 10pt}&
&\smash{\vrule height 10pt}&
&\smash{\vrule height 10pt}&
&\smash{\vrule height 10pt}&
&\smash{\vrule height 10pt}&
&\smash{\vrule height 10pt}&
&\smash{\vrule height 10pt}&
&\smash{\vrule height 10pt}&
&\smash{\vrule height 10pt}&
&\smash{\vrule height 10pt}&
\cr
\hline
height 0.2em depth 0.2em\gobble&&
&\vrule depth0pt&
&\vrule depth0pt&
&\vrule depth0pt&
&\vrule depth0pt&
&\vrule depth0pt&
&\vrule depth0pt&
&\vrule depth0pt&
&\vrule depth0pt&
&\vrule depth0pt&
&\vrule depth0pt&
&\vrule depth0pt&
&\vrule depth0pt&
&\vrule depth0pt&
\cr
height 0pt\gobble&&
&\DSp\putdown{10}\DSp&
&\TSp\putdown{108}\TSp&
&\TSp\putdown{365}\TSp&
&\TSp\putdown{606}\TSp&
&\TSp\putdown{535}\TSp&
&\TSp\putdown{254}\TSp&
&\DSp\putdown{61}\DSp&
&\QSp\putdown{6}\QSp&
&&
&&
&&
&&
&&\cr
&-9&
\group{0}{0}{10}{}&\rightarrowfill&
\group{0}{0}{118}{}&\rightarrowfill&
\group{1}{0}{474}{\bf}&\rightarrowfill&
\group{1}{1}{972}{\bf}&\rightarrowfill&
\group{0}{0}{1141}{}&\rightarrowfill&
\group{0}{0}{789}{}&\rightarrowfill&
\group{0}{0}{315}{}&\rightarrowfill&
\group{0}{0}{67}{}&\rightarrowfill&
\group{0}{0}{6}{}&\smash{\vrule height 15pt}&
&\smash{\vrule height 15pt}&
&\smash{\vrule height 15pt}&
&\smash{\vrule height 15pt}&
&\smash{\vrule height 15pt}&
\cr
height 0.3em\gobble&&
&\smash{\vrule height 10pt}&
&\smash{\vrule height 10pt}&
&\smash{\vrule height 10pt}&
&\smash{\vrule height 10pt}&
&\smash{\vrule height 10pt}&
&\smash{\vrule height 10pt}&
&\smash{\vrule height 10pt}&
&\smash{\vrule height 10pt}&
&\smash{\vrule height 10pt}&
&\smash{\vrule height 10pt}&
&\smash{\vrule height 10pt}&
&\smash{\vrule height 10pt}&
&\smash{\vrule height 10pt}&
\cr
\hline
height 0.2em depth 0.2em\gobble&&
&\vrule depth0pt&
&\vrule depth0pt&
&\vrule depth0pt&
&\vrule depth0pt&
&\vrule depth0pt&
&\vrule depth0pt&
&\vrule depth0pt&
&\vrule depth0pt&
&\vrule depth0pt&
&\vrule depth0pt&
&\vrule depth0pt&
&\vrule depth0pt&
&\vrule depth0pt&
\cr
height 0pt\gobble&&
&\QSp\putdown{5}\QSp&
&\DSp\putdown{38}\DSp&
&\DSp\putdown{80}\DSp&
&\DSp\putdown{88}\DSp&
&\DSp\putdown{48}\DSp&
&\DSp\putdown{12}\DSp&
&\QSp\putdown{1}\QSp&
&&
&&
&&
&&
&&
&&\cr
&-11&
\group{0}{0}{5}{}&\rightarrowfill&
\group{0}{0}{43}{}&\rightarrowfill&
\group{1}{1}{119}{\bf}&\rightarrowfill&
\group{0}{0}{168}{}&\rightarrowfill&
\group{0}{0}{136}{}&\rightarrowfill&
\group{0}{0}{60}{}&\rightarrowfill&
\group{0}{0}{13}{}&\rightarrowfill&
\group{0}{0}{1}{}&\smash{\vrule height 15pt}&
&\smash{\vrule height 15pt}&
&\smash{\vrule height 15pt}&
&\smash{\vrule height 15pt}&
&\smash{\vrule height 15pt}&
&\smash{\vrule height 15pt}&
\cr
height 0.3em\gobble&&
&\smash{\vrule height 10pt}&
&\smash{\vrule height 10pt}&
&\smash{\vrule height 10pt}&
&\smash{\vrule height 10pt}&
&\smash{\vrule height 10pt}&
&\smash{\vrule height 10pt}&
&\smash{\vrule height 10pt}&
&\smash{\vrule height 10pt}&
&\smash{\vrule height 10pt}&
&\smash{\vrule height 10pt}&
&\smash{\vrule height 10pt}&
&\smash{\vrule height 10pt}&
&\smash{\vrule height 10pt}&
\cr
\hline
height 0.2em depth 0.2em\gobble&&
&\vrule depth0pt&
&\vrule depth0pt&
&\vrule depth0pt&
&\vrule depth0pt&
&\vrule depth0pt&
&\vrule depth0pt&
&\vrule depth0pt&
&\vrule depth0pt&
&\vrule depth0pt&
&\vrule depth0pt&
&\vrule depth0pt&
&\vrule depth0pt&
&\vrule depth0pt&
\cr
height 0pt\gobble&&
&\QSp\putdown{1}\QSp&
&\QSp\putdown{4}\QSp&
&\QSp\putdown{6}\QSp&
&\QSp\putdown{4}\QSp&
&\QSp\putdown{1}\QSp&
&&
&&
&&
&&
&&
&&
&&
&&\cr
&-13&
\group{0}{0}{1}{}&\rightarrowfill&
\group{1}{0}{6}{\bf}&\rightarrowfill&
\group{0}{0}{10}{}&\rightarrowfill&
\group{0}{0}{10}{}&\rightarrowfill&
\group{0}{0}{5}{}&\rightarrowfill&
\group{0}{0}{1}{}&\smash{\vrule height 15pt}&
&\smash{\vrule height 15pt}&
&\smash{\vrule height 15pt}&
&\smash{\vrule height 15pt}&
&\smash{\vrule height 15pt}&
&\smash{\vrule height 15pt}&
&\smash{\vrule height 15pt}&
&\smash{\vrule height 15pt}&
\cr
height 0.3em\gobble&&
&\smash{\vrule height 10pt}&
&\smash{\vrule height 10pt}&
&\smash{\vrule height 10pt}&
&\smash{\vrule height 10pt}&
&\smash{\vrule height 10pt}&
&\smash{\vrule height 10pt}&
&\smash{\vrule height 10pt}&
&\smash{\vrule height 10pt}&
&\smash{\vrule height 10pt}&
&\smash{\vrule height 10pt}&
&\smash{\vrule height 10pt}&
&\smash{\vrule height 10pt}&
&\smash{\vrule height 10pt}&
\cr
\hline
\dblhline
}}$$

 \end{footnotesize}
\caption{Ranks of $\mc{H}^{i,j}$ and $\mc{C}^{i,j}$ and ranks of the chain differentials
for either the Kinoshita-Terasaka knot or the Conway knot (both have the same Khovanov
homology)}\label{t3}
\end {sidewaystable}

\end{document}